\documentclass[11pt]{amsart}
\usepackage[a4paper,margin=27mm]{geometry}
\usepackage{amsmath,amssymb,amsthm,mathtools,booktabs,microtype,enumitem}
\usepackage{pdfsync}

\usepackage{mathrsfs,graphicx,latexsym,tikz,color,euscript}
\usepackage{amsfonts,amscd}
\usepackage{epstopdf}
\usepackage{upgreek}
\usepackage{comment}
\usepackage{tikz-cd}
\usepackage[colorlinks=true,linkcolor=blue,citecolor=blue,urlcolor=blue]{hyperref}

\newtheorem{theorem}{Theorem}[section]
\newtheorem{proposition}[theorem]{Proposition}
\newtheorem{lemma}[theorem]{Lemma}
\newtheorem{corollary}[theorem]{Corollary}
\theoremstyle{remark}
\newtheorem{remark}[theorem]{Remark}

\newcommand{\R}{\mathbb R}
\newcommand{\C}{\mathbb C}
\newcommand{\Z}{\mathbb Z}
\newcommand{\dd}{\,\mathrm d}
\newcommand{\norm}[1]{\left\lVert#1\right\rVert }
\newcommand{\ip}[2]{  \langle \langle #1,#2 \rangle\rangle  }
\newcommand{\U}{U_\alpha}
\newcommand{\Jh}{\mathcal J_h}

\newcommand{\al}{\alpha}

\newcommand{\ep}{\varepsilon}

\newcommand{\dl}{\delta}

\title[hyperbolic scattering with strong force potentials]{Hyperbolic scattering in the planar $n$-body problem with homogeneous strong force potentials}

\author{Junhao Lie}
\address{Chern Institute of Mathematics and LPMC, Nankai University, Tianjin, China}
\email{liejh@mail.nankai.edu.cn}

\author{Guowei Yu}
\address{Chern Institute of Mathematics and LPMC, Nankai University, Tianjin, China}
\email{yugw@nankai.edu.cn}

\thanks{This work is supported by NSFC (No. 12671224), Nankai Zhide Foundation and the Fundamental Research Funds for the Central Universities. }

\begin{document}

\begin{abstract}
For the planar $n$-body problem with homogeneous strong force potentials, given an arbitrary pair of normalized collision-free configurations $(a, b)$, we prove that for any positive energy $h$, there are infinitely many $h$-energy scattering orbits between $a$ and $b$.   
\end{abstract}

\maketitle

\section{Introduction and main results} \label{sec;intro}

Consider the \(n\)-body problem, with $n\ge2$, in $\R^2$ with the following homogeneous potential function of degree $-\al$ 
\[
 \U(q)=\sum_{1\le i<j\le n}\frac{m_im_j}{|q_i-q_j|^\alpha}, \;\; \al >0
\]
where $m_i>0$ and $q_i \in \R^2$ represent the mass and location of the $i$-th body, for each $1 \le i \le n$. We say $U_{\al}$ is a strong force potential when $\al \ge 2$, and is a weak force potential when $0 < \al < 2$. The classical Newtonian potential corresponds to $\al =1$. 

Let $\langle \cdot , \cdot \rangle$ be the Euclidean inner product. We may introduce the following mass-weighted inner product on $(\R^2)^n$
\[
 \ip{x}{y}=\sum_{i=1}^n m_i \langle x_i, y_i \rangle.
\] 
The induced mass-weighted norm will be denoted by $\| \cdot \|$. 

The motion of the masses is described by the following equation 
\begin{equation}\label{eq:system}
 \ddot q=\nabla_M\U(q),
\end{equation}
where \(\nabla_M\) is the gradient with respect to the mass-weighted inner product. 

It is well known that the linear momentum is a first integral of the $n$-body problem. Hence there is no loss of generality in fixing the center of mass at the origin by setting
\[
 E=\left\{q\in(\R^2)^n:\sum_{i=1}^n m_iq_i=0\right\}.
\]
Then \(E^*=E\setminus\Delta\) is the set of collision-free configurations, where \(\Delta\) is the set of collision configurations 
\[
 \Delta=\bigcup_{i<j}\{q\in E:q_i=q_j\}.
\]

Since the $n$-body problem is a time-independent Hamiltonian system, the following energy is conserved along any solution of \eqref{eq:system} 
\[
 h=\tfrac12\norm{\dot q}^2-\U(q).
\]
We are interested in the final motions of the system. In particular those have the following asymptotic behaviors, when time goes to positive or negative infinity. 
\[
 |q_i(t)-q_j(t)|\to\infty,\qquad
 |\dot q_i(t)-\dot q_j(t)|\to |v_i-v_j|>0
 \quad(i\ne j).
\]
 Following Chazy's terminology \cite{Chazy22}, they will be called \emph{hyperbolic motions}. For hyperbolic motions in more general potentials, see \cite{FKM21}. 

To give a more precise description of hyperbolic solutions, denote the unit mass sphere and its collision-free subset by
$$  S=\{s\in E:\norm{s}=1\}, \qquad  S^*=S\setminus\Delta. $$
For any positive energy $h$, if $q|_{(-\infty, t_0]}$, resp. $q|_{[t_0, +\infty)}$, is an $h$-energy hyperbolic solution, there exists a collision-free normalized configuration $a \in S^*$, resp. $b \in S^*$, such that the following asymptotic expansions hold
\begin{equation} \label{eq;asym}
q(t) = -\sqrt{2h} t a + o(|t|),  \quad t \to -\infty \quad \text{resp.} \quad  q(t) = \sqrt{2h} t b + o(t), \quad t \to +\infty.
\end{equation}
Depending on the value of $\al$, more precise asymptotic expansions can be obtained. Here we are mainly interested in the normalized limiting configurations and do not pursue these refinements. 

\begin{theorem}
\label{thm;fix-initial-data}
Let $\alpha>0$, $h>0$, and $q_0\in E$, with $q_0\in E^*$ if $\alpha\ge2$. For any $a\in S^*$ resp. $b\in S^*$,  there exist a collision-free $h$-energy hyperbolic solution $q|_{(-\infty, t_0]}$, resp. $q|_{[t_0, +\infty)}$, satisfying $q(t_0) = q_0$ the corresponding asymptotic expansion in \eqref{eq;asym}. 
\end{theorem}
Results of this type were proved in \cite{MV20,LYZ21,PT24,Yu24} for different ranges of $\alpha$, using variational methods to obtain the solutions as action minimizers. In particular, Theorem 1.4 and Corollary 1.5 of \cite{LYZ21} cover all $\alpha>0$, including strong force potentials. An important ingredient for weak force potentials is Marchal's lemma (see \cite{C02,FT04}), which excludes collisions in the interior of the parameter interval of fixed-end minimizers. For strong force potentials, every finite-action curve on a compact time interval is collision-free, including at its endpoints. This is why we require $q_0\in E^*$ when $\alpha\ge2$.

Meanwhile instead of considering the problem with a given initial configuration, for an arbitrary pair of normalized configurations $(a, b) \in S^* \times S^*$, we may ask if there is an $h$-energy solution $q|_{\R}$ satisfying both asymptotic expansions in \eqref{eq;asym}. Such a solution, if it exists, will be called \emph{bi-hyperbolic}. As we can see, this is deeply connected with the scattering theory of particles, and if such an $h$-energy bi-hyperbolic solution exists, we say there is an $h$-energy scattering between $a$ and $b$. 

For the Newtonian potential and more general weak force potentials, the prescribed scattering problem remains largely open. Relevant results include \cite{DMMY20,Moeckel26}. For further discussion of scattering in the Newtonian $n$-body problem, see Montgomery's book \cite{Mont25}.

For the simplified $n$-center problem, where $n$ bodies are fixed and only one is free to move, there are many more results on scattering; see \cite{Knauf92,BDT17,BDP18,BBD18}. For a relativistic two-center problem, see \cite{SC26}.

In this paper we prove the following result for strong force potentials.
\begin{theorem}\label{thm:main}
When \(\alpha\ge2\), for any positive $h$ and normalized collision-free configurations $(a, b) \in S^* \times S^*$, there exist infinitely many $h$-energy scattering orbits between $a$ and $b$, distinct modulo time translation. 
\end{theorem}

Our approach is again variational. Unlike the fixed-initial-configuration problem, the variational study of the scattering problem is more difficult due to the lack of coercivity, or more precisely the possibility of the escape of the minimizing sequence to infinity. A possible way to avoid this is to consider minimization problems under suitable topological constraints. For the planar $n$-body problem, such topological constraints can be described by the braid group. 

This has been used by Montgomery \cite{Mont98} in the variational study of relative periodic orbits for the planar $n$-body problem with strong force potentials. Using Gordon's notion of tied classes, Montgomery imposed \emph{the tied condition}, which forces the length of a curve to go to infinity, if the curve escapes to infinity. This provides the needed control for minimizing sequences of closed curves in a fixed class, whose lengths have a finite upper bound. 

However, for the problem of scattering, the lengths of all admissible paths are infinite. Hence the tied condition, which is qualitative in nature, is not enough to force coercivity. A quantitative condition is needed, and we will explain what it is now.  

Given any two points \(a,b\in S^*\), let \(\Gamma(a, b)\) denote a relative homotopy class of rectifiable paths on \(S^*\) connecting $a$ and $b$ with endpoints fixed during homotopies. Define the infimum of the length of all these paths as
\begin{equation}\label{eq:lengthgap}
 \ell(\Gamma(a, b))=\inf_{s\in\Gamma(a, b)}\mathcal L(s),\qquad
 \mathcal L(s)=\int\norm{s'(u)}\dd u.
\end{equation}
The coercive condition we need is 
\begin{equation} 
\ell(\Gamma(a, b)) > \pi. \label{eq:coer}
\end{equation}

\begin{theorem} \label{thm:class}
 When $\al \ge 2$, for any $h >0$ and $(a, b) \in S^* \times S^*$, if a relative homotopy class $\Gamma(a, b)$ satisfies condition \eqref{eq:coer}, then there exists an $h$-energy bi-hyperbolic solution $q|_{\R}$ between $a$ and $b$ with the normalized curve $q/\|q\|$, extended continuously at both ends, belonging to $\Gamma(a, b)$. 
 \end{theorem} 

 Theorem \ref{thm:main} will follow directly from the above result, once we show that there are infinitely many different homotopy classes $\Gamma(a, b)$ satisfying condition \eqref{eq:coer}, which is precisely the result of Proposition \ref{prop:coercive} below. To explain it, define the rotation path 
$$ R_a(u)=e^{2\pi iu}a, \quad u\in[0,1]. $$
Since $R_a$ is one full rotation of all bodies, its homotopy class $[R_a]$ is the positive full twist. 
$R^{-1}_a$ will represent the inverse path of $R_a$. Fix an arbitrary smooth path $\eta$ in $S^*$ between $a$ and $b$. For each $k \in \Z$, define a homotopy class  
$$
 \Gamma_k(a, b)=[R_a^k*\eta], 
$$
where $R_a^k$ denotes the $k$-fold concatenation for $k>0$, the constant loop for $k=0$, and the $|k|$-fold concatenation of $R_a^{-1}$ for $k<0$, and $*$ denotes chronological path concatenation. 

\begin{proposition}\label{prop:coercive}
For any \(a,b\in S^*\), there exists an integer \(k_0\), depending only on the path \(\eta\) and the masses, such that for all \(|k|\ge k_0\), \(\ell(\Gamma_k(a, b))>\pi\).  
\end{proposition}




\section{Positive energy solutions and geodesics of the Jacobi metric }

Given any positive energy constant $h$, define the following Jacobi metric on $E^*$
\[g_h=2(h+\U)g, \; \; g \text{ the mass metric.} \] 
We use absolutely continuous parametrizations of rectifiable curves throughout. For such a curve \(\xi\) in \(E^*\), define its Jacobi length as 
\begin{equation}\label{eq:JM}
 \Jh(\xi)=\int f(\xi)\norm{\xi'}\dd u,\qquad
 f(q)=\sqrt{2(h+\U(q))}.
\end{equation}
By the Jacobi--Maupertuis principle, nonconstant geodesics on \((E^*,g_h)\), up to a reparametrization, are $h$-energy solutions of \eqref{eq:system}. 

To see this, let $\xi$ be a nonconstant geodesic with a regular parameter $u$. Recall that the Euler--Lagrange equation of the functional \(\Jh\) is
\[
 \frac{\dd}{\dd u}\left(f(\xi)\frac{\xi'}{\norm{\xi'}}\right)=\norm{\xi'}\nabla_M f(\xi)
\]
where $'$ is the derivative with respect to $u$. Now introduce the new time parameter $t$ by 
\begin{equation} \label{eq;t}
 \dd t/\dd u=\norm{\xi'}/f(\xi)>0
\end{equation}
and set \(q(t):=\xi(u(t))\), where \(u(t)\) is the inverse reparametrization.  Then
\[
 \dot q=\frac{\dd u}{\dd t}\xi'=f(\xi)\frac{\xi'}{\norm{\xi'}}
\]
where a dot denotes differentiation with respect to \(t\). Differentiating  again with respect to \(t\) shows that $q(t)$ is a solution of Newton's equation
\[ \ddot q=\frac{f}{\norm{\xi'}}\frac{\dd}{\dd u}\!\left(f(\xi)\frac{\xi'}{\norm{\xi'}}\right)=f\nabla_M f=\nabla_M\U. \]
Accordingly, the parameter $t$ defined above will be referred to as the \emph{physical time} of the geodesic. 
Moreover the following energy identity holds 
\[
 \tfrac12\norm{\dot q}^2-\U(q)=\tfrac12 f(\xi)^2-\U(q)=h. 
\]

Finally, in this parametrization,
\begin{equation}\label{eq:lengthtime}
 \dd\Jh=f(\xi)\norm{\xi'}\dd u=f(\xi)^2\dd t=2(h+\U(q))\dd t\ge2h\dd t.
\end{equation}

Although $E^*$ is an open subset of $E$, the strong force potential implies \((E^*,g_h)\) is a smooth complete Riemannian manifold, as shown by the next lemma. 

\begin{lemma}\label{prop:complete}
When \(h>0\) and \(\alpha\ge 2\), \((E^*,g_h)\) is a complete smooth
Riemannian manifold.
\end{lemma}

\begin{proof}
Since \(U_\alpha\in C^\infty(E^*)\) and \(h+U_\alpha>0\) on \(E^*\), the
Riemannian metric
\[
g_h=2(h+U_\alpha)g
\]
is smooth on \(E^*\). Hence it remains to prove metric completeness. Let
\(d_h\) be the distance induced by the Jacobi length
\[
\Jh(\xi)=\int \sqrt{2(h+U_\alpha(\xi))}\,\|\xi'\|\,du.
\]

We first establish a pairwise-distance estimate. Fix \(i<j\) and write
\[
d_{ij}(q):=|q_i-q_j|,\qquad q = (q_i)_{i=1}^n \in E^*.
\]
Let \(\xi:[a,b]\to E^*\) be an absolutely continuous curve and put
\[
d_{ij}(u):=d_{ij}(\xi(u)).
\]
Since \(\xi(u)\in E^*\), we have \(d_{ij}(u)>0\). At every \(u\) at which
\(\xi'(u)\) exists,
\[
d'_{ij}(u)
=
\frac{\bigl\langle \xi_i(u)-\xi_j(u),\xi_i'(u)-\xi_j'(u)\bigr\rangle}
{|\xi_i(u)-\xi_j(u)|},
\]
and therefore, by the Cauchy--Schwarz inequality in \(\mathbb R^2\),
\[
|d'_{ij}(u)|\le |\xi_i'(u)-\xi_j'(u)|.
\]
On the other hand,
\[
\begin{aligned}
|\xi_i'(u)-\xi_j'(u)|
&\le |\xi_i'(u)|+|\xi_j'(u)|\\
&\le
\left(\frac1{m_i}+\frac1{m_j}\right)^{1/2}
\left(m_i|\xi_i'(u)|^2+m_j|\xi_j'(u)|^2\right)^{1/2}\\
&\le
\left(\frac1{m_i}+\frac1{m_j}\right)^{1/2}\|\xi'(u)\|.
\end{aligned}
\]
Consequently,
\[
\|\xi'(u)\|
\ge
\left(\frac1{m_i}+\frac1{m_j}\right)^{-1/2}|d'_{ij}(u)|
\]
for almost every \(u\). Moreover,
\[
\sqrt{2(h+U_\alpha(\xi(u)))}
\ge
\sqrt{2m_i m_j}\,d_{ij}(u)^{-\alpha/2}.
\]
Hence
\begin{equation}\label{eq:curve-pairwise}
\begin{aligned}
\Jh(\xi)
&\ge
\frac{\sqrt{2m_i m_j}}
{\left(m_i^{-1}+m_j^{-1}\right)^{1/2}}
\int_a^b d_{ij}(u)^{-\alpha/2}|d'_{ij}(u)|\,du.
\end{aligned}
\end{equation}
Set
 \begin{equation}
 \label{eq;G-a} G_\alpha(r)=
\begin{cases}
\log r, & \alpha=2,\\[2mm]
\dfrac{r^{1-\alpha/2}}{1-\alpha/2}, & \alpha>2.
\end{cases}
 \end{equation}
Then \(G_\alpha'(r)=r^{-\alpha/2}\), and \(G_\alpha(r)\to-\infty\) as
\(r\searrow0\). Therefore
\[
\int_a^b d_{ij}(u)^{-\alpha/2}|d'_{ij}(u)|\,du
\ge
\left|
G_\alpha(d_{ij}(\xi(b)))-G_\alpha(d_{ij}(\xi(a)))
\right|.
\]
Taking the infimum over all curves in \(E^*\) joining two points
\(p,q\in E^*\), we obtain
\begin{equation}\label{eq:pairwise-distance}
d_h(p,q)\ge
\frac{\sqrt{2m_i m_j}}
{\left(m_i^{-1}+m_j^{-1}\right)^{1/2}}
\left|
G_\alpha(d_{ij}(q))-G_\alpha(d_{ij}(p))
\right|.
\end{equation}

Now let \(\{x_k\}_{k\ge1}\subset E^*\) be a Cauchy sequence with respect to
\(d_h\). Since \(h+U_\alpha\ge h\), for every curve \(\eta\) joining two
points \(p,q\in E^*\) we have
\[
\Jh(\eta)
\ge
\sqrt{2h}\int\|\eta'\|\,du
\ge
\sqrt{2h}\,\|p-q\|.
\]
Thus
\[
d_h(p,q)\ge \sqrt{2h}\,\|p-q\|.
\]
Hence \(\{x_k\}\) is Cauchy, and therefore bounded, in the mass norm
\(\|\cdot\|\). Since \(E\) is finite dimensional, \(\{x_k\}\) is contained
in a bounded ball of \(E\).

Next we show that \(\{x_k\}\) cannot approach the collision set. Fix
\(i<j\). By \eqref{eq:pairwise-distance}, the scalar sequence
\[
\left\{G_\alpha(d_{ij}(x_k))\right\}_{k\ge1}
\]
is Cauchy in \(\mathbb R\), hence converges to a finite number. If
\[
\liminf_{k\to\infty}d_{ij}(x_k)=0,
\]
then along a subsequence \(d_{ij}(x_k)\to0\), which would imply
\[
G_\alpha(d_{ij}(x_k))\longrightarrow-\infty,
\]
a contradiction. Therefore
\[
\liminf_{k\to\infty}d_{ij}(x_k)>0.
\]
Since there are only finitely many pairs \(i<j\), there exist
\(\delta>0\) and \(k_0\) such that
\[
d_{ij}(x_k)\ge \delta
\qquad
\text{for all }i<j\text{ and all }k\ge k_0.
\]

Because \(\{x_k\}\) is bounded in the finite-dimensional space \(E\), there
is a subsequence \(x_{k_l}\) converging in the mass norm to some
\(x_\infty\in E\). Discarding finitely many terms if necessary,  we may assume
\(k_l\ge k_0\) for all \(l\). Letting \(l\to\infty\) in the separation
estimate gives
\[
d_{ij}(x_\infty)\ge \delta>0
\qquad
\text{for all }i<j.
\]
Hence \(x_\infty\in E^*\).

It remains to prove that this subsequence converges to \(x_\infty\) in
\(d_h\). Since \(x_\infty\in E^*\) and \(\Delta\) is closed in \(E\), there
exists \(\rho>0\) such that the closed mass-norm ball
\[
\overline B_\rho(x_\infty)
=
\{q\in E:\|q-x_\infty\|\le\rho\}
\]
is contained in \(E^*\). On this compact ball the function
\[
f(q)=\sqrt{2(h+U_\alpha(q))}
\]
is continuous and positive. Set
\[
\Lambda:=\max_{q\in \overline B_\rho(x_\infty)} f(q)<\infty.
\]
If \(p,q\in B_{\rho/2}(x_\infty)\), then the straight segment joining
\(p\) and \(q\) lies in \(\overline B_\rho(x_\infty)\subset E^*\), and its
Jacobi length is at most \(\Lambda\|p-q\|\). Hence
\[
d_h(p,q)\le \Lambda\|p-q\|.
\]
For all sufficiently large \(l\), we have
\(x_{k_l}\in B_{\rho/2}(x_\infty)\), and therefore
\[
d_h(x_{k_l},x_\infty)
\le
\Lambda\|x_{k_l}-x_\infty\|\longrightarrow0.
\]
Thus a subsequence of \(\{x_k\}\) converges to \(x_\infty\) in \(d_h\).

Since \(\{x_k\}\) is \(d_h\)-Cauchy and one of its subsequences converges
to \(x_\infty\) in \(d_h\), the whole sequence converges to \(x_\infty\) in
\(d_h\). Indeed, given \(\varepsilon>0\), choose \(N\) such that
\[
d_h(x_k,x_l)<\frac{\varepsilon}{2}
\qquad
\text{for all }k,l\ge N.
\]
Choose \(k_{l_0}\ge N\) with
\[
d_h(x_{k_{l_0}},x_\infty)<\frac{\varepsilon}{2}.
\]
Then for every \(k\ge N\),
\[
d_h(x_k,x_\infty)
\le
d_h(x_k,x_{k_{l_0}})
+
d_h(x_{k_{l_0}},x_\infty)
<\varepsilon.
\]
Therefore every \(d_h\)-Cauchy sequence converges in \(E^*\). Hence
\((E^*,d_h)\) is complete. Together with the smoothness of \(g_h\), this
proves that \((E^*,g_h)\) is a complete smooth Riemannian manifold.
\end{proof}

\section{Some estimates of the Jacobi length}

For the rest of the paper, assume $\alpha\ge2$. In this section, fix a pair of collision-free normalized configurations \(a,b\in S^*\), and a relative homotopy class \(\Gamma(a, b)\) of paths in \(S^*\) connecting \(a\) and \(b\) with \(\ell(\Gamma(a, b))>\pi\), and a piecewise smooth representative $\gamma\in\Gamma(a,b)$. For simplicity, write $\Gamma=\Gamma(a,b)$. Such representatives exist by piecewise smooth approximation in the open manifold $S^*$.

For $R>1$, let $\Gamma_R$ be the set of rectifiable curves in $E^*$ from $Ra$ to $Rb$ with the corresponding normalized curves belonging to $\Gamma$. Since $(E^*, g_h)$ is a complete smooth Riemmannian manifold by Lemma \ref{prop:complete}, the infimum of $\Jh$ in $\Gamma_R$ is a minimum and we denote such a minimizer by $q_R$ for the rest of the section. As a geodesic, $q_R$ can be naturally parameterized by its arc length. Meanwhile we can also introduce a physical-time $t$ as given in \eqref{eq;t} to parameterize $q_R$.   


\begin{lemma} \label{lem:exterior}
Let \(\xi\) be an absolutely continuous curve whose two endpoint norms equal \(R\), with minimum norm \(\rho>0\) and  \(\mathcal{L}(\xi/\|\xi\|) =L\). Then
\begin{equation}\label{eq:exterior}
 \mathcal L(\xi)\ge2\left(\sqrt{R^2-\rho^2}-\rho\arccos\frac\rho R\right)+\rho L
 \ge2R+\rho(L-\pi).
\end{equation}
\end{lemma}
\begin{proof}
Write \(\xi=rs\), where \(r=\norm{\xi}\) and \(\norm s=1\). By a direct computation
\[
 \norm{\xi'}^2=r'^2+r^2\norm{s'}^2.
\]
Since \(r\ge\rho\), the following is a unit vector
\[ \bigl(\sqrt{1-\rho^2/r^2},\rho/r\bigr).\] 
Applying the Cauchy--Schwarz inequality to \((|r'|,r\norm{s'})\) and the above unit vector 
\[
 \norm{\xi'}=\sqrt{r'^2+r^2\norm{s'}^2}
 \ge\sqrt{1-\frac{\rho^2}{r^2}}\,|r'|+\frac{\rho}{r}\,r\norm{s'}
 =\sqrt{1-\frac{\rho^2}{r^2}}\,|r'|+\rho\norm{s'}.
\]
Let $F_{\rho}(r)$ be a function of $r$ with domain $[\rho, \infty)$, 
\[
 F_\rho(r)=\sqrt{r^2-\rho^2}-\rho\arccos\frac\rho r.
\]
Then \(F_\rho(\rho)=0\) and by direct differentiation 
\[F_\rho'(r)=\sqrt{1-\rho^2/r^2}. \]

Write the parameter interval of \(\xi\) as \([\mu_-,\mu_+]\), and choose
\(u_0\in[\mu_-,\mu_+]\) such that
\[
r(u_0)=\rho=\min_{u\in[\mu_-,\mu_+]}r(u).
\]
At almost every \(u\), the chain rule gives
\[
(F_\rho\circ r)'(u)
=
F_\rho'(r(u))r'(u)
=
\sqrt{1-\frac{\rho^2}{r(u)^2}}\,r'(u).
\]
As a result, 
\[
\|\xi'(u)\|
\ge
\left|(F_\rho\circ r)'(u)\right|
+\rho\|s'(u)\|.
\]
Integrating this inequality over \([\mu_-,\mu_+]\), we obtain
\[
\begin{aligned}
\mathcal L(\xi)
&=
\int_{\mu_-}^{\mu_+}\|\xi'(u)\|\,du\\
&\ge
\int_{\mu_-}^{\mu_+}\left|(F_\rho\circ r)'(u)\right|\,du
+
\rho\int_{\mu_-}^{\mu_+}\|s'(u)\|\,du\\
&=
\int_{\mu_-}^{u_0}\left|(F_\rho\circ r)'(u)\right|\,du
+
\int_{u_0}^{\mu_+}\left|(F_\rho\circ r)'(u)\right|\,du
+\rho L.
\end{aligned}
\]
Since \(r(\mu_-)=R\) and \(r(u_0)=\rho\),  
\[
\int_{\mu_-}^{u_0}\left|(F_\rho\circ r)'(u)\right|\,du
\ge
\left|F_\rho(r(u_0))-F_\rho(r(\mu_-))\right|
=
F_\rho(R)-F_\rho(\rho)
=
F_\rho(R),
 \]
Similarly, since \(r(u_0)=\rho\) and \(r(\mu_+)=R\),
\[
\int_{u_0}^{\mu_+}\left|(F_\rho\circ r)'(u)\right|\,du
\ge
F_\rho(R)-F_\rho(\rho)
=
F_\rho(R).
\]
Therefore
\[
\mathcal L(\xi)
\ge
2F_\rho(R)+\rho L
=
2\left(
\sqrt{R^2-\rho^2}
-\rho\arccos\frac{\rho}{R}
\right)
+\rho L.
\]
This proves the first inequality.

For the second inequality, it suffices to prove that
\begin{equation}\label{eq:lemma31-second}
F_\rho(R)-R+\frac{\pi\rho}{2}\ge0
\qquad(R\ge\rho).
\end{equation}
Notice that the function \(F_\rho(R)-R+\frac{\pi\rho}{2}\) is nonincreasing w.r.t.
\(R\), as its derivative is
\begin{equation}\label{eq:lemma31-derivative}
\sqrt{1-\frac{\rho^2}{R^2}}-1\le0,
\end{equation}
Moreover,
\begin{equation}\label{eq:lemma31-limit-one}
F_\rho(R)-R
=
\sqrt{R^2-\rho^2}-R
-\rho\arccos\frac{\rho}{R}
\longrightarrow
-\frac{\pi\rho}{2}
\qquad(R\to\infty),
\end{equation}
thus
\begin{equation}\label{eq:lemma31-limit-two}
\lim_{R\to\infty}
\left(
F_\rho(R)-R+\frac{\pi\rho}{2}
\right)
=0.
\end{equation}
Combining this with the monotonicity in \eqref{eq:lemma31-derivative} proves the required inequality \eqref{eq:lemma31-second}.

\end{proof}

\begin{lemma}\label{lem:Jh-bounds}
There is a constant \(C>0\), independent of \(R\), such that
\begin{equation}
\label{eq:upper}
 2\sqrt{2h}\,R+\sqrt{2h}\,\rho_R(\ell(\Gamma)-\pi)\le\Jh(q_R)\le2\sqrt{2h}\,R+C,
\end{equation}
where \(\rho_R=\min\norm{q_R}\).
\end{lemma}
\begin{proof}
For the upper bound, consider the curve from \(Ra\) to \(Rb\) obtained by concatenating the incoming path
\[
 \rho_a(r)=(R+1-r)a,\qquad r\in[1,R],
\]
which goes from \(\rho_a(1)=Ra\) to \(\rho_a(R)=a\), the representative \(\gamma\), and the outgoing path
\[
 \rho_b(r)=rb,\qquad r\in[1,R],
\]
which goes from \(\rho_b(1)=b\) to \(\rho_b(R)=Rb\). By the minimizing property of \(q_R\),
\[
 \Jh(q_R)\le\Jh(\rho_a)+\Jh(\gamma)+\Jh(\rho_b).
\]

To estimate \(\Jh(\rho_b)\), notice that \(\norm{\rho_b'(r)}=\norm b=1\), and homogeneity of the potential gives \(\U(\rho_b(r))=r^{-\alpha}\U(b)\). As a result,
\[
 0\le\sqrt{2(h+r^{-\alpha}\U(b))}-\sqrt{2h}
 =\frac{2r^{-\alpha}\U(b)}{\sqrt{2(h+r^{-\alpha}\U(b))}+\sqrt{2h}}
 \le\frac{r^{-\alpha}\U(b)}{\sqrt{2h}},
\]
and hence
\[
 \Jh(\rho_b)=\int_1^R\sqrt{2\bigl(h+r^{-\alpha}\U(b)\bigr)}\dd r
 \le\sqrt{2h}\,(R-1)+\frac{\U(b)}{\sqrt{2h}}\int_1^\infty r^{-\alpha}\dd r.
\]
Similarly, by substituting \(a\) for \(b\), we obtain
\[
 \Jh(\rho_a)\le\sqrt{2h}\,(R-1)+\frac{\U(a)}{\sqrt{2h}}\int_1^\infty r^{-\alpha}\dd r.
\]
Since $\Jh(\gamma)$ is finite, there is a constant \(C>0\), independent of \(R\), such that
\begin{equation}
 \Jh(q_R)\le2\sqrt{2h}\,R+C.
\end{equation}

For the lower bound, recall that by \eqref{eq:JM}, \(f\ge\sqrt{2h}\). Hence 
$$\Jh(q_R)\ge\sqrt{2h}\,\mathcal L(q_R). 
$$
Let \(L_R\) be the spherical length of the normalized curve \(q_R/\norm{q_R}\). Since the normalized curve is a representative of \(\Gamma\), \(L_R\ge\ell(\Gamma)>\pi\). Applying Lemma~\ref{lem:exterior} to \(q_R\) with \(\rho=\rho_R >0\) gives
\[
 \mathcal L(q_R)\ge2R+\rho_R(L_R-\pi)\ge2R+\rho_R(\ell(\Gamma)-\pi).
\] 
Together with $\Jh(q_R)\ge\sqrt{2h}\,\mathcal L(q_R)$ and the upper bound just proved, this implies
\[
 2\sqrt{2h}\,R+\sqrt{2h}\,\rho_R(\ell(\Gamma)-\pi)\le\Jh(q_R)\le2\sqrt{2h}\,R+C.
\]
\end{proof}

The next proposition gives the desired coercivity in our proof of the main theorem. 
\begin{proposition}\label{prop:minradius}
There is \(B>0\), independent of \(R\), such that \(\rho_R=\min\norm{q_R}\le B\).
\end{proposition}
\begin{proof}
By Lemma~\ref{lem:Jh-bounds},
\[
 2\sqrt{2h}\,R+\sqrt{2h}\,\rho_R(\ell(\Gamma)-\pi)\le\Jh(q_R)\le2\sqrt{2h}\,R+C.
\]
Since \(\ell(\Gamma)>\pi\), subtracting \(2\sqrt{2h}\,R\) and dividing by \(\sqrt{2h}\,(\ell(\Gamma)-\pi)>0\) gives
\begin{equation}
\label{eq;B} \rho_R\le\frac{C}{\sqrt{2h}\,(\ell(\Gamma)-\pi)}=:B.
\end{equation}
\end{proof}

\section{Proof of Theorem \ref{thm:class}}

We retain the notation of the previous section. In particular let $B$ and $C$ be the constants given in \eqref{eq;B} and \eqref{eq:upper}. For each $R>\max\{1,B\}$, let $q_R$ be the minimizer obtained in the previous section. After a suitable translation of the time parameter, we may assume the domain of physical-time parameter $t$ of $q_R$ is $[-T_R^-, T_R^+]$ and   
$$ \rho_R = \|q_R(0)\| = \min \{ \|q_R(t)\|: t \in [-T_R^-, T_R^+] \}.  $$

First we prove a lemma which controls the normalized configurations near the endpoints of the minimizers. For this we set 
\[
 D=\frac{C}{\sqrt{2h}}+2B.
\]

\begin{lemma} \label{lem:cone}
Let \(q_R(t)=r(t)s(t)\) with $r = \|q_R\|$ and $s = q_R/\|q_R\|$. Then  
\begin{equation}\label{eq:cone}
 r(1-\ip a s)\le D,\qquad \norm{s-a}^2\le\frac{2D}{r}, \; \forall t \in [-T_R^-, 0].
\end{equation}
The same inequalities hold for $t \in [0, T_R^+]$ after replacing $a$ with $b$.
\end{lemma}
\begin{remark}
The geometric argument below uses only the length bound \eqref{eq:upper} and the uniform bound on the minimum radius. It therefore applies beyond strong force potentials whenever these two bounds are available, and can be seen as a generalization and simplification of Lemmas 2.1 and 3.1 in \cite{Yu24}. 
\end{remark}

\begin{proof}
Fix an arbitrary $t \in [-T_R^-, 0]$. Recall that  \(\norm{q_R(0)}\le B\) by Proposition~\ref{prop:minradius}. Then 
\begin{equation} \label{eq:three-piece-length}
\begin{aligned}
\mathcal L(q_R) & \ge \|Ra-q_R(t)\|+\|q_R(t)-q_R(0)\|+\|q_R(0)-Rb\| \\
& \ge \|Ra-q_R(t)\|+(r-B)+(R-B),
\end{aligned}
\end{equation}
Meanwhile \(\Jh(q_R)\le2\sqrt{2h}\,R+C\) by \eqref{eq:upper}, and \(\Jh(q_R)\ge\sqrt{2h}\,\mathcal L(q_R)\), so
\begin{equation}\label{eq:length-control}
 \mathcal L(q_R)\le\frac{\Jh(q_R)}{\sqrt{2h}}\le2R+\frac{C}{\sqrt{2h}}.
\end{equation}
Combining \eqref{eq:three-piece-length} and \eqref{eq:length-control},
\begin{equation}\label{eq:length-combined}
 2R+\frac{C}{\sqrt{2h}}\ge\|Ra-q_R(t)\|+(r-B)+(R-B).
\end{equation}

Since \(\norm a=1\), by the Cauchy--Schwarz inequality
\[ \ip{a}{q_R(t)}^2\le\norm{q_R(t)}^2=r^2\]
As a result, 
\begin{equation}\label{eq:radial-bound}
 \begin{aligned}
 \|Ra-q_R(t)\|^2 &  =R^2-2R\ip{a}{q_R(t)}+\norm{q_R(t)}^2 \\
 & \ge R^2-2R\ip{a}{q_R(t)}+ \ip{a}{q_R(t)}^2 \\
 &=\bigl(R-\ip{a}{q_R(t)}\bigr)^2,
 \end{aligned}
\end{equation}
which implies \(\|Ra-q_R(t)\| \ge R-\ip{a}{q_R(t)}\). Substituting this lower bound into \eqref{eq:length-combined} yields
\[
 r-\ip{a}{q_R(t)}\le\frac{C}{\sqrt{2h}}+2B=D.
\]
This proves the first inequality in \eqref{eq:cone}. Using it and the following identity gives us the second inequality
\[\norm{s-a}^2=2(1-\ip a s).\]

The result for $t \in [0, T_R^+]$ follows similarly. 
\end{proof}

A direct corollary of the above result is 

\begin{corollary}
There exist positive constants $K>\max\{1,B\}$, $C_1, C_2$ and $\dl_0$  independent of the choice of $q_R$, such that the following inequalities hold for all $q_R(t)$ when $\|q_R(t)\| \ge K$, 
\begin{equation}\label{eq:tailbounds}
 \U(q_R)\le C_1 \|q_R\|^{-\alpha},\qquad
 \norm{\nabla_M\U(q_R)}\le C_2 \|q_R\|^{-\alpha-1},\qquad
 |q_{R,i}-q_{R,j}|\ge \dl_0 \|q_R\|.
\end{equation}
\end{corollary}

For each $q_R(t)$, let $I_R(t)=\|q_R(t)\|^2$ be its moment of inertia. By the Lagrange--Jacobi identity 
\begin{equation}\label{eq:virial}
 \ddot{I}_R=4h-2(\alpha-2)\U(q_R).
\end{equation}
Further assume that $K$ is large enough such that 
$$ 2(\alpha-2)C_1 K^{-\alpha}\le2h$$
where $C_1$ is as given in \eqref{eq:tailbounds}, and this implies 
\begin{equation}
\label{eq;ddot-I} \ddot{I}_R(t) \ge 2h>0, \; \text{ when } \|q_R(t)\| \ge K. 
\end{equation}
From now on, let $R>K$. Recall that $\|q_R(0)\|=\rho_R<K$. By \eqref{eq;ddot-I}, $\{t:\|q_R(t)\|\le K\}$ is a single closed interval $[\tau_R^-,\tau_R^+]$ containing $0$, where
\[
 \tau_R^- = \max\{t\in[-T_R^-,0]:\|q_R(t)\|=K\},\qquad
 \tau_R^+ = \min\{t\in[0,T_R^+]:\|q_R(t)\|=K\}.
\]
Indeed, a bounded exterior excursion would have $I_R>K^2$ inside and $I_R=K^2$ at both endpoints, contradicting strict convexity. Moreover, the radial derivatives at the outgoing and incoming crossings have nonnegative and nonpositive signs, respectively. Integrating \eqref{eq;ddot-I} gives
 \begin{equation}
 \label{eq:growth} \begin{aligned}
  \norm{q_R(t)}^2&\ge K^2+h(t-\tau_R^+)^2 &&(\tau_R^+\le t\le T_R^+),\\
  \norm{q_R(t)}^2&\ge K^2+h(t-\tau_R^-)^2 &&(-T_R^-\le t\le\tau_R^-). 
 \end{aligned}
 \end{equation}

The next lemma shows that both endpoints of the physical-time interval tend to infinity in absolute value as $R\to\infty$. 
\begin{lemma}\label{lem:core}
For every \(R>K\),
\[
\tau_R^+-\tau_R^-\le \frac{C+2\sqrt{2h}K}{2h},
\qquad
T_R^\pm\ge \frac{R-K}{\sqrt{2(h+C_1K^{-\alpha})}},
\]
where \(C\) and \(C_1\) are the constants in \eqref{eq:upper} and
\eqref{eq:tailbounds}, respectively.
\end{lemma}

\begin{proof}
Recall that for the physical-time \(t\), 
\[
\dd\Jh=2(h+U_\alpha(q_R(t)))\,\dd t.
\]
By the endpoint condition for \(q_R\) and the definition of \(\tau_R^\pm\),
\[
\|q_R(-T_R^-)\|=\|q_R(T_R^+)\|=R,
\qquad
\|q_R(\tau_R^\pm)\|=K.
\]
Moreover by \eqref{eq:growth}, 
\[
\|q_R(t)\|\ge K
\qquad
t\in[-T_R^-,\tau_R^-]\cup[\tau_R^+,T_R^+].
\]

We first estimate \(\tau_R^+-\tau_R^-\). Since $\|\dot q_R(t)\|\ge\sqrt{2h}$
and $\dd\Jh=\|\dot q_R(t)\|^2\,\dd t$
on the interval \([-T_R^-,\tau_R^-]\), we have
\[
\Jh(q_R|_{[-T_R^-,\tau_R^-]})
\ge
\sqrt{2h}\int_{-T_R^-}^{\tau_R^-}\|\dot q_R(t)\|\,dt
\ge
\sqrt{2h}(R-K).
\]
Similarly,
\[
\Jh(q_R|_{[\tau_R^+,T_R^+]})
\ge
\sqrt{2h}(R-K).
\]
By \eqref{eq:upper},
\[
\Jh(q_R)\le 2\sqrt{2h}R+C.
\]
Hence
\[
\begin{aligned}
\Jh(q_R|_{[\tau_R^-,\tau_R^+]})
&=
\Jh(q_R)
-
\Jh(q_R|_{[-T_R^-,\tau_R^-]})
-
\Jh(q_R|_{[\tau_R^+,T_R^+]})\\
&\le
2\sqrt{2h}R+C
-
2\sqrt{2h}(R-K)\\
&=
C+2\sqrt{2h}K.
\end{aligned}
\]
On the other hand, by \eqref{eq:lengthtime},
\[
\Jh(q_R|_{[\tau_R^-,\tau_R^+]})
\ge
2h(\tau_R^+-\tau_R^-).
\]
Therefore
\[
\tau_R^+-\tau_R^-\le \frac{C+2\sqrt{2h}K}{2h}.
\]

Next we estimate \(T_R^\pm\). On both intervals
\[
[-T_R^-,\tau_R^-],
\qquad
[\tau_R^+,T_R^+],
\]
we have \(\|q_R(t)\|\ge K\). Hence, by \eqref{eq:tailbounds},
\[
U_\alpha(q_R(t))\le C_1K^{-\alpha}.
\]
The energy identity gives
\[
\|\dot q_R(t)\|^2
=
2(h+U_\alpha(q_R(t)))
\le
2(h+C_1K^{-\alpha}),
\]
so
\[
\|\dot q_R(t)\|
\le
\sqrt{2(h+C_1K^{-\alpha})}.
\]
On \([-T_R^-,\tau_R^-]\),
\[
R-K
\le
\int_{-T_R^-}^{\tau_R^-}\|\dot q_R(t)\|\,dt
\le
\sqrt{2(h+C_1K^{-\alpha})}\,(T_R^-+\tau_R^-).
\]
Since \(\tau_R^-\le0\),
\[
T_R^-
\ge
T_R^-+\tau_R^-
\ge
\frac{R-K}{\sqrt{2(h+C_1K^{-\alpha})}}.
\]
Similarly, on \([\tau_R^+,T_R^+]\),
\[
R-K
\le
\int_{\tau_R^+}^{T_R^+}\|\dot q_R(t)\|\,dt
\le
\sqrt{2(h+C_1K^{-\alpha})}\,(T_R^+-\tau_R^+),
\]
and since \(\tau_R^+\ge0\),
\[
T_R^+
\ge
T_R^+-\tau_R^+
\ge
\frac{R-K}{\sqrt{2(h+C_1K^{-\alpha})}}.
\]
This proves the required estimates.
\end{proof}

Now we can give a proof of Theorem~\ref{thm:class}. 
\begin{proof}[Proof of Theorem~\ref{thm:class}]
Recall that under condition \eqref{eq:coer}, Proposition~\ref{prop:minradius} implies
\[
\|q_R(0)\|\le B
\]
for every \(R>K\).

We first show that \(q_R(0)\) stays uniformly away from the collision set.
From the proof of Lemma~\ref{lem:core} we have
\[
\Jh(q_R|_{[\tau_R^-,\tau_R^+]})\le C+2\sqrt{2h}K.
\]
Since \(0\in[\tau_R^-,\tau_R^+]\), the same bound holds for the subarc
\(q_R|_{[\tau_R^-,0]}\). At \(t=\tau_R^-\), inequality \eqref{eq:tailbounds} gives
\begin{equation}
\label{eq;min-dist} |q_{R,i}(\tau_R^-)-q_{R,j}(\tau_R^-)|\ge \delta_0K
\qquad(i<j).
\end{equation}
Fix \(i<j\). Applying the estimate \eqref{eq:curve-pairwise} to the curve
\(q_R|_{[\tau_R^-,0]}\), and using the uniform upper bound for its Jacobi
length, we obtain a constant \(\eta_{ij}>0\), independent of \(R\), such
that
\[
|q_{R,i}(0)-q_{R,j}(0)|\ge \eta_{ij}
\qquad(R>K).
\]
Taking
\[
\eta=\min_{i<j}\eta_{ij}>0,
\]
we conclude that
\[
\|q_R(0)\|\le B,
\qquad
|q_{R,i}(0)-q_{R,j}(0)|\ge\eta
\quad(i<j)
\]
for every \(R>K\). Thus all \(q_R(0)\) lie in a fixed compact subset of
\(E^*\). The energy identity
\[
\|\dot q_R(0)\|^2=2(h+U_\alpha(q_R(0)))
\]
therefore gives a constant, independent of \(R\), which bounds
\(\|\dot q_R(0)\|\).

Passing to a subsequence \(R_k\to\infty\), we may assume
\[
q_{R_k}(0)\longrightarrow q_0,
\qquad
\dot q_{R_k}(0)\longrightarrow p_0.
\]
Since the above compact subset is contained in \(E^*\), we have
\(q_0\in E^*\).

Let \(q\) be the solution of \eqref{eq:system} with 
\[
q(0)=q_0,
\qquad
\dot q(0)=p_0.
\]
By continuity, the energy of $q$ must be $h$. 

Next we claim that there exists $\varepsilon_*>0$, independent of $R$ and $t$, such that 
\begin{equation} \label{eq;qR-collision-free}
|q_{R,i}(t)-q_{R,j}(t)|\ge\varepsilon_*, \quad t\in[-T_R^-,T_R^+],\quad R>K,\quad 1\le i<j\le n. 
\end{equation}
Let $d_{ij}(t)=|q_{R,i}(t)-q_{R,j}(t)|$.
For $t\in[\tau_R^-,\tau_R^+]$, the subarc $q_R|_{[\tau_R^-,t]}$ has Jacobi length at most $ C+2\sqrt{2h}K$. Hence \eqref{eq:curve-pairwise} and \eqref{eq;min-dist} give
\[
\begin{aligned}
 G_\alpha(d_{ij}(t)) & \ge G_\alpha(d_{ij}(\tau_R^-))-
 \frac{\sqrt{m_i^{-1}+m_j^{-1}}}{\sqrt{2m_im_j}}\bigl(C+2\sqrt{2h}K\bigr) \\
 &\ge G_\alpha(\delta_0K)-
 \frac{\sqrt{m_i^{-1}+m_j^{-1}}}{\sqrt{2m_im_j}}\bigl(C+2\sqrt{2h}K\bigr).
\end{aligned}
\]
Since $G_\alpha$ is increasing and $G_\alpha(r)\to-\infty$ as $r\searrow0$, there exists $\ep_{ij}>0$, independent of $R$, such that  
$$ d_{ij}(t)\ge\varepsilon_{ij}, \quad t \in [\tau_R^-, \tau_R^+]. $$
Meanwhile  by \eqref{eq:tailbounds} 
$$ d_{ij}(t)\ge\delta_0K, \quad t \in [-T_R^-, \tau_R^-] \cup [\tau_R^+, T_R^+].$$
Hence \eqref{eq;qR-collision-free} follows with $\varepsilon_*:=\min\{\delta_0K,\min_{i<j}\varepsilon_{ij}\}>0$.

By \eqref{eq;qR-collision-free} and the energy identity,
\[
 \|\dot q_R(t)\|^2\le
 2\left(h+\sum_{i<j}m_im_j\varepsilon_*^{-\alpha}\right)=:V^2,
 \qquad \|q_R(t)\|\le B+V|t|.
\]
By the above result, at any finite time $t$, $(q,\dot q)(t)$ would remain in a compact subset of $E^*\times E$, so the smooth differential equation would extend the solution. Thus $q$ is collision-free and defined for all $t\in\R$.

For every fixed $T>0$, Lemma~\ref{lem:core} gives $T_{R_k}^\pm\ge T$ for all sufficiently large $k$. Continuous dependence for \eqref{eq:system} now gives uniform convergence of $q_{R_k}$ and its first two derivatives to $q$ and its corresponding derivatives on $[-T,T]$.

Next we show $q(t)$ satisfies the prescribed asymptotic behavior in both time directions. Since
\(0\in[\tau_{R_k}^-,\tau_{R_k}^+]\), Lemma~\ref{lem:core} implies that
\(\tau_{R_k}^\pm\) are bounded. Passing to a further subsequence, we may
assume
\[
\tau_{R_k}^-\longrightarrow\tau^-,
\qquad
\tau_{R_k}^+\longrightarrow\tau^+.
\]
Passing to the limit in \eqref{eq:growth}, we obtain
\[
\|q(t)\|^2\ge K^2+h(t-\tau^+)^2
\qquad(t\ge\tau^+),
\]
and
\[
\|q(t)\|^2\ge K^2+h(t-\tau^-)^2
\qquad(t\le\tau^-).
\]
In particular,
\[
\|q(t)\|\longrightarrow\infty
\qquad(|t|\to\infty).
\]

Set $s(t)=\frac{q(t)}{\|q(t)\|}$. For fixed \(t\le0\), Lemma~\ref{lem:cone} and the convergence on \([-T,T]\) give
\[
\|q(t)\|\bigl(1-\langle\langle a,s(t)\rangle\rangle\bigr)\le D,
\qquad
\|s(t)-a\|^2\le\frac{2D}{\|q(t)\|}.
\]
Letting \(t\to-\infty\) and using \(\|q(t)\|\to\infty\), we get
\[
s(t)\longrightarrow a
\qquad(t\to-\infty).
\]
Similarly, the corresponding estimates with \(b\) give
\[
s(t)\longrightarrow b
\qquad(t\to+\infty).
\]

We now verify hyperbolicity. The growth estimates and the collision-free limiting configurations imply that $|q_i(t)-q_j(t)|\ge c|t|$ for all $i<j$ and sufficiently large $|t|$, with some $c>0$. Consequently,
\[
 U_\alpha(q(t))=O(|t|^{-\alpha}),\qquad
 \|\ddot q(t)\|=O(|t|^{-\alpha-1}).
\]
Integration gives $\dot q(t)=v_\pm+O(|t|^{-\alpha})$ as $t\to\pm\infty$. The energy identity yields $\|v_\pm\|=\sqrt{2h}$, while $q(t)/t\to v_\pm$ and the limits of $s(t)$ identify
\[
 v_-=-\sqrt{2h}\,a,\qquad v_+=\sqrt{2h}\,b.
\]
Thus $q$ is bi-hyperbolic and satisfies \eqref{eq;asym}.  

By the above result, we can extend the normalized curve $s$ continuously to the endpoints $a$ and $b$. Moreover as for each $R_k$, the normalized curve $q_{R_k}/\|q_{R_k}\|$ belongs to the class $\Gamma_{a, b}$, by the above result $q_{R_k}/\|q_{R_k}\|$ converges uniformly to $s$ in $S^*$ with the endpoints fixed, so $s$ must belong to $\Gamma_{a, b}$ as well. 



\end{proof}


\section{Proof of Proposition \ref{prop:coercive}}

Fix the reference path $\eta:a\to b$ from Section~\ref{sec;intro}, and write $\Gamma_k=\Gamma_k(a,b)=[R_a^k*\eta]$. We shall prove the stronger statement
\[
 \ell(\Gamma_k)\longrightarrow\infty\qquad\text{as }|k|\longrightarrow\infty.
\]
The main point is that bounded-length paths can approach partial collisions, but the braiding confined near a partial collision does not contribute to the invariant used below.

First assume $n=2$. In this case, the centered configuration space $E$ is a real two-dimensional space with the only collision configuration at $0$. Thus $S^*=S$ is a unit circle in the mass metric. Choose a lift $\theta_\eta:[0,1]\to\R$ of the angular coordinate along $\eta$. Every path $s\in\Gamma_k$ has angular increment $\theta_\eta(1)-\theta_\eta(0)+2\pi k$, and hence
\[
 \mathcal L(s)\ge
 \left|\theta_\eta(1)-\theta_\eta(0)+2\pi k\right|.
\]
This proves the proposition for $n=2$.

For the rest of this section assume $n\ge3$, and identify the plane with $\C$. The ordered configuration space is
\[
 \operatorname{Conf}_n(\C)=\{q\in\C^n:q_i\ne q_j\text{ for }i\ne j\}.
\]
For $q\in\operatorname{Conf}_n(\C)$, let $c(q)=\sum_i m_iq_i/\sum_i m_i$ and $\mathbf1=(1,\ldots,1)$. Translation and radial normalization give
\begin{equation}\label{eq:confprod}
 \begin{aligned}
 \operatorname{Conf}_n(\C)&\cong\C\times(0,\infty)\times S^*,\\
 q&\longmapsto\left(c(q),\|q-c(q)\mathbf1\|,
 \frac{q-c(q)\mathbf1}{\|q-c(q)\mathbf1\|}\right).
 \end{aligned}
\end{equation}
In particular $E^*\cong(0,\infty)\times S^*$. Since the first two factors in \eqref{eq:confprod} are contractible, $\pi_1(S^*)$ identifies with the pure braid group $P_n$, viewed as a subgroup of the braid group $B_n$; see \cite{Birman74}. Under this identification, $[R_a]$ is the positive full twist $Z$.

Recall that a quasimorphism on a group $G$ is a function $f:G\to\R$ whose defect
\[
 D(f):=\sup_{g,h\in G}|f(gh)-f(g)-f(h)|
\]
is finite. It is homogeneous if $f(g^k)=kf(g)$ for every $g\in G$ and $k\in\Z$.

\begin{theorem}[Fractional Dehn twist coefficient; \cite{Malyutin,Feller}]\label{thm:fdtc}
For $n\ge3$, there exists a homogeneous quasimorphism $\omega:B_n\to\R$ such that
\begin{equation}\label{eq:fdtc}
 \begin{gathered}
 |\omega(xy)-\omega(x)-\omega(y)|\le1,
 \qquad \omega(x^k)=k\omega(x),\\
 \omega(Z)=1,\qquad \omega|_{B_{n-1}}=0,
 \end{gathered}
\end{equation}
where $B_{n-1}$ is the standard subgroup supported on $n-1$ strands.
\end{theorem}

Malyutin introduced the twist number in \cite{Malyutin}; Feller's characterization in \cite{Feller} gives uniqueness under the conditions above. We use only existence and the properties in \eqref{eq:fdtc}.

\begin{lemma}\label{lem:qm}
The quasimorphism $\omega$ satisfies
\begin{equation}\label{eq;k01}
 \omega(1)=0,\qquad \omega(g^{-1})=-\omega(g)\quad(g\in B_n).
\end{equation}
It is invariant under conjugation and additive on commuting elements. Moreover, for $N\ge1$,
\begin{equation}\label{eq:defect}
 \left|\omega(g_1\cdots g_N)-\sum_{j=1}^N\omega(g_j)\right|\le N-1.
\end{equation}
\end{lemma}
\begin{proof}
The identities in \eqref{eq;k01} follow from homogeneity with exponents $0$ and $-1$. Applying the defect bound twice gives, for every positive integer $k$,
\[
 |\omega(ga^kg^{-1})-\omega(g)-\omega(a^k)-\omega(g^{-1})|\le2.
\]
The terms involving $g$ cancel. Since $ga^kg^{-1}=(gag^{-1})^k$, division by $k$ and passage to the limit prove conjugacy invariance.

If $ab=ba$, then $(ab)^k=a^kb^k$. The defect bound applied to $a^k,b^k$, followed by homogeneity, gives
\[
 k|\omega(ab)-\omega(a)-\omega(b)|\le1.
\]
Letting $k\to\infty$ proves additivity. Repeated application gives additivity for any finite collection of pairwise commuting elements. Finally, induction using the two-factor defect bound proves \eqref{eq:defect}.
\end{proof}

Fix the identification $\pi_1(S^*,a)=P_n$. A loop based at another point is transported to $a$ by a connecting path. Different choices conjugate the resulting element, so Lemma~\ref{lem:qm} makes its $\omega$-value independent of this choice. This convention will be used for all local loops below.

\begin{lemma}\label{lem:local}
For every $p\in S$, every sufficiently small open geodesic ball $V$ centered at $p$ has the property that $V\cap S^*$ is dense in $V$ and any two of its points can be joined in $V\cap S^*$ by a piecewise smooth path.
\end{lemma}
\begin{proof}
Let \(d_M\) denote the intrinsic distance on \(S\) induced by the mass metric.
For any \(p\in S\), choose \(r>0\) sufficiently small, in particular
\(r<\pi/2\), and let
\[
V=\{s\in S:d_M(s,p)<r\}
\]
be the open geodesic ball in \(S\) centered at \(p\) with respect to the
mass metric. Choose \(r\) small enough that \(V\) is geodesically convex;
that is, for any two points of \(V\), the minimizing geodesic segment
joining them in \(S\) lies in \(V\).

Let \(x,y\in V\cap S^*\). We first find a curve in \(E\) joining \(x\)
and \(y\) that does not pass through any collision point. For \(i<j\),
set
\[
L_{ij}=\{q\in E:q_i=q_j\}.
\]
Then \(L_{ij}\) is a linear subspace of \(E\) of codimension two. Denote by
\[
\operatorname{aff}(L_{ij}\cup\{x\})
\]
the smallest affine subspace of \(E\) containing \(L_{ij}\) and \(x\).
Since \(L_{ij}\) is a linear subspace, this is precisely the linear span
of \(L_{ij}\) and \(x\). Because \(x\in S^*\), we have
\(x\notin L_{ij}\), and consequently
\(\operatorname{aff}(L_{ij}\cup\{x\})\) has codimension at least one
in \(E\).

Therefore the finite union
\[
\bigcup_{i<j}
\bigl(\operatorname{aff}(L_{ij}\cup\{x\})\cup
\operatorname{aff}(L_{ij}\cup\{y\})\bigr)
\]
is a finite union of proper affine subspaces of \(E\), so its complement
\[
E\setminus\bigcup_{i<j}
\bigl(\operatorname{aff}(L_{ij}\cup\{x\})\cup
\operatorname{aff}(L_{ij}\cup\{y\})\bigr)
\]
is dense in \(E\). We may choose a point \(z\) in this complement
sufficiently close to \(p\).

If the line segment \([x,z]\) met some \(L_{ij}\), then, since
\(x\notin L_{ij}\), it would meet \(L_{ij}\) at an interior point.
Solving for \(z\) from such an intersection would give
\(z\in\operatorname{aff}(L_{ij}\cup\{x\})\), contradicting the choice
of \(z\). Hence \([x,z]\) meets no \(L_{ij}\). The same argument shows that
\([z,y]\) meets no \(L_{ij}\). Thus the polygonal segemnt
\([x,z]\cup [z,y]\) is collision-free.

Since \(z\) is sufficiently close to \(p\), we have
\(z/\|z\|\in V\). Applying the radial projection
\(q\mapsto q/\|q\|\) to the line segments \([x,z]\) and \([z,y]\), we
obtain two geodesic segments in \(S\), joining \(x\) to \(z/\|z\|\) and
\(z/\|z\|\) to \(y\), respectively. Because \(L_{ij}\) is homogeneous,
the radial projection of a collision-free curve is again collision-free.
Hence both geodesic segments lie in \(V\cap S^*\). Therefore any two
points \(x,y\in V\cap S^*\) can be joined in \(V\cap S^*\) by a piecewise
smooth path.
\end{proof}

\begin{lemma}\label{lem:vanish}
Every $p\in S$ has an open geodesic-ball neighborhood $V_p$ such that $V_p\cap S^*$ is dense and path connected as in Lemma~\ref{lem:local}, and every loop in $V_p\cap S^*$ has $\omega$-value zero. These properties remain valid when the ball is shrunk about $p$.
\end{lemma}
\begin{proof}
Partition the labels into clusters according to equality of positions at $p$:
\[
 \{1,\ldots,n\}=I_1\sqcup\cdots\sqcup I_r,
 \qquad p_i=p_j\iff i,j\in I_\nu\text{ for some }\nu.
\]
There are at least two clusters. Otherwise centering would give $p=0$, contradicting $\|p\|=1$. In particular, each cluster contains at most $n-1$ labels.

Choose disjoint open planar disks $D_1,\ldots,D_r$ around the distinct cluster positions. By coordinate continuity, a sufficiently small ball $V_p$ satisfies
\[
 s\in V_p,\ i\in I_\nu\quad\Longrightarrow\quad s_i\in D_\nu.
\]
Shrink it further so that Lemma~\ref{lem:local} applies.

Let $\gamma$ be a loop in $V_p\cap S^*$, based at $x$. Regard it as a loop in the full ordered configuration space. For each cluster $I_\nu$, define $\widehat\gamma_\nu$ by letting its labels follow their motions in $\gamma$ and holding all other labels fixed at their positions in $x$. These loops need not be centered or normalized, which is why we work in the full configuration space here.

The subspace in which, for every $\nu$, the labels of $I_\nu$ remain in $D_\nu$ is a product of the corresponding ordered configuration spaces of the disks. Its fundamental group is a product as well. Consequently, after inclusion in the full braid group,
\[
 [\gamma]=[\widehat\gamma_1]\cdots[\widehat\gamma_r],\qquad
 [\widehat\gamma_\nu][\widehat\gamma_\mu]
 =[\widehat\gamma_\mu][\widehat\gamma_\nu].
\]
Each factor is supported in a disk containing at most $n-1$ strands, and hence belongs to a conjugate of a standard braid subgroup contained in $B_{n-1}$. Equations \eqref{eq:fdtc} and conjugacy invariance give $\omega([\widehat\gamma_\nu])=0$. Additivity on commuting elements now yields $\omega([\gamma])=0$. The same disk containment holds for every smaller ball, proving the final assertion.
\end{proof}

\begin{proposition}\label{prop:topology}
For every finite $L\ge0$, the numbers $\omega([s*\bar\eta])$ are bounded as $s:a\to b$ ranges over paths in $S^*$ of spherical length at most $L$. Here $\bar\eta$ denotes the reversed path. 

\end{proposition}
\begin{proof}
Suppose the assertion fails. Choose paths $s_\nu:a\to b$ with $\mathcal L(s_\nu)\le L$ and
\[
 |\omega([s_\nu*\bar\eta])|\longrightarrow\infty.
\]
Parametrize them with constant speed on $[0,1]$, allowing constant paths when their length is zero. Their Lipschitz constants are at most $L$. By compactness of $S$ and the Arzel\`a--Ascoli theorem, a subsequence converges uniformly to a path $s_\infty:[0,1]\to S$. The limit may meet the collision set.

Choose a subdivision $0=u_0<\cdots<u_N=1$ and balls $V_j$ supplied by Lemma~\ref{lem:vanish} such that
\[
 s_\infty([u_{j-1},u_j])\subset V_j\quad(1\le j\le N).
\]
Such a subdivision exists by a Lebesgue number for the inverse-image cover of $[0,1]$. Each displayed image is compact inside its open ball. Uniform convergence therefore implies the same containment for $s_\nu$ when $\nu$ is sufficiently large.

For each interior vertex, choose a small ball $O_j$ centered at $s_\infty(u_j)$ whose closure is contained in $V_j\cap V_{j+1}$, and to which Lemma~\ref{lem:local} applies. Choose a fixed point $y_j\in O_j\cap S^*$. For sufficiently large $\nu$, both $y_j$ and $s_\nu(u_j)$ lie in this path-connected set, so choose a connector
\[
 \rho_{\nu,j}:y_j\longrightarrow s_\nu(u_j)
 \quad\text{inside }O_j\cap S^*.
\]
Set $y_0=a$, $y_N=b$, and take the endpoint connectors to be constant. Because $y_{j-1},y_j\in V_j\cap S^*$, choose fixed paths $\zeta_j:y_{j-1}\to y_j$ there, and put $\zeta=\zeta_1*\cdots*\zeta_N$.

Let $s_{\nu,j}$ be the restriction of $s_\nu$ to $[u_{j-1},u_j]$ and define
\[
 \beta_{\nu,j}=\rho_{\nu,j-1}*s_{\nu,j}*\bar\rho_{\nu,j},
 \qquad \lambda_{\nu,j}=\beta_{\nu,j}*\bar\zeta_j.
\]
The first path joins $y_{j-1}$ to $y_j$; the second is a loop at $y_{j-1}$ contained in $V_j\cap S^*$. Hence $\omega([\lambda_{\nu,j}])=0$ by Lemma~\ref{lem:vanish}.

Write $P_0$ for the constant path at $a$ and $P_{j-1}=\zeta_1*\cdots*\zeta_{j-1}$. Cancellation of each connector with its reverse, and then of the reference paths, gives the identity in $\pi_1(S^*,a)$
\[
 [s_\nu*\bar\zeta]
 =\prod_{j=1}^N[P_{j-1}*\lambda_{\nu,j}*\bar P_{j-1}].
\]
Each factor has $\omega$-value zero by conjugacy invariance. The defect estimate \eqref{eq:defect} gives
\[
 |\omega([s_\nu*\bar\zeta])|\le N-1.
\]
Finally, $[s_\nu*\bar\eta]=[s_\nu*\bar\zeta][\zeta*\bar\eta]$, so one more application of the defect bound yields
\[
 |\omega([s_\nu*\bar\eta])|
 \le N+|\omega([\zeta*\bar\eta])|.
\]
The right-hand side is independent of $\nu$, contradicting the chosen sequence.

\end{proof}

Now we explain how the above result can be used to prove 
\begin{equation}\label{eq:ellinfty}
 \ell(\Gamma_k)\longrightarrow\infty\qquad\text{as }|k|\longrightarrow\infty.
\end{equation}
Notice that $[s*\bar\eta]=Z^k$ whenever $s\in\Gamma_k$, and therefore $\omega([s*\bar\eta])=k$. If \eqref{eq:ellinfty} failed, there would be a constant $L_0$ and integers $k_\nu$ with $|k_\nu|\to\infty$ and $\ell(\Gamma_{k_\nu})\le L_0$. By the definition of the infimum, choose $s_\nu\in\Gamma_{k_\nu}$ with $\mathcal L(s_\nu)\le L_0+1$. This contradicts the boundedness result given by the above proposition, and this finishes our proof of Proposition~\ref{prop:coercive}.

The classes $\Gamma_k$ are distinct: for $n\ge3$ their loop invariants equal $k$, and for $n=2$ their winding numbers differ by $k$. Applying Theorem~\ref{thm:class} to every sufficiently large positive or negative $k$ therefore gives infinitely many scattering orbits. Time translation preserves the compactified normalized path class, so these orbits are distinct modulo time translation. This completes the proof of Theorem~\ref{thm:main}.

\medskip
\noindent\textbf{Declaration of AI usage.} AI tools assisted in developing the results in the last section. The authors take full responsibility for their correctness.  

\bibliographystyle{abbrv}
\bibliography{RefScattering}

\end{document}